\documentclass[11pt, a4paper]{article}
\usepackage{geometry}
\usepackage{amsmath, amssymb, amscd, amsthm, amsfonts}
\usepackage{bm, mathrsfs,dsfont}
\usepackage{tikz}
\usepackage{color}
\usepackage{cases}
\usepackage{float}  
\usepackage{bbding}

\newtheorem{theorem}{Theorem}[section]

\title{\textbf{Isolated singularities of flat metrics on Riemann surfaces}}
\author{Jin Li and Bin Xu}
\date{}
\begin{document}
\maketitle

\paragraph{Abstract}
Robert Bryant ({\small Theorie des varietes
minimales et applications, 1988, 154: 321-347}) proved that an isolated singularity of a conformal metric of positive constant curvature on a Riemann surface is a conical one. Using Complex Analysis, we find all of the local models for an isolated singularity of a flat metric whose area satisfies some polynomial growth condition near the singularity. In particular, we show that an isolated singularity of a flat metric with finite area is also a conical one.



\section{Introduction}
\paragraph{}
Let $\Sigma$ be a Riemann surface and $p$ a point on $\Sigma$. A conformal metric $\mathrm{d}{\sigma}^{2}$ on $\Sigma$ has a \emph{conical singularity} at $p$ with \emph{cone angle} $2\pi\alpha>0$ if in a neighborhood of $p$, $\mathrm{d}{\sigma}^{2}=e^{2\varphi}|\mathrm{d}z|^2$, where $z$ is a local complex coordinate defined in the neighborhood of $p$ with $z(p)=0$ and $\varphi-(\alpha-1)\mathrm{ln}|z|$ is continuous in the neighborhood.

Since the singularity $p$ is isolated, by choosing a suitable complex coordinate centered at $p$, we can assume that $\mathrm{d}{\sigma}^{2}$ is a conformal metric on the punctured disk ${\Delta}^{*}=\{\omega\in\mathbb{C}\mid0<|\omega|<1\}$.

Bryant proved that an isolated singularity of a conformal metric of Gauss curvature $+1$ with finite area must be a conical singularity, and gave an explicit expression of the metric near it. See the following theorem.

\begin{theorem}
Let $\mathrm{d}{\sigma}^{2}$ be a conformal metric of Gauss curvature $+1$ on the punctured disk ${\Delta}^{*}=\{\omega\in\mathbb{C}\mid0<|\omega|<1\}$. Suppose moreover that the $\mathrm{d}{\sigma}^{2}$-area of ${\Delta}^{*}$ is finite. Then there exists a local holomorphic coordinate $z$ on  $U_{\varepsilon}=\{\omega\in\mathbb{C}\mid |\omega|<\varepsilon \}$ for some $\varepsilon>0$ with $z(0)=0$, and a real number $\beta>-1$ so that, on $U_{\varepsilon}$, we have
\begin{equation}
\mathrm{d}{\sigma}^{2} \Big|_{U_\varepsilon} =\frac{4(\beta+1)^2  |z|^{2\beta}}{(1+|z|^{2(\beta+1)})^2} |\mathrm{d}z|^2.  \notag
\end{equation}
Moreover, $\beta$ is unique, and $z$ is unique up to a replacement by $\lambda z$ where $|\lambda|=1$.
\end{theorem}

We omit the proof, in which Bryant used elementary value-distribution theory. We want to show in the following theorem that in the flat case, under a more relaxing but still controllable area condition, the possible isolated singularities include, but are not limited to conical ones.

\begin{theorem}{\rm (Main Theorem)}   Let $\mathrm{d}{\sigma}^{2}$ be a conformal metric of Gauss curvature zero on  the punctured disk ${\Delta}^{*}=\left\{\omega\in\mathbb{C}\mid0<|\omega|<1\right\}$. Suppose moreover that the $\mathrm{d}{\sigma}^{2}$-area near the origin $\omega=0$ satisfies polynomial growth condition. More precisely, there exist $0<R<1$ and $M>0, N\geq0,$ which are independent of $R$ such that
\begin{equation}
\textrm{Area}\bigg(\Delta\Big(0, \frac{1}{r}, R\Big)\bigg)\le M{r}^{N}, \forall {r} > \frac{1}{R},    \notag
\end{equation}
where $\Delta\big(0, \frac{1}{r}, R\big)=\left\{ \omega\in\mathbb{C}\mid\frac{1}{r}<|\omega|<R \right\}$ denotes the annulus with inner radius $\frac{1}{r}$ and outer radius $R$. Then, there exists a holomorphic coordinate transformation $\omega\mapsto{z}(\omega)$ on  $U_{\varepsilon}=\{\omega\in\mathbb{C}\mid |\omega|<\varepsilon \}$ for some $\varepsilon>0$ with $z(0)=0$ so that on $U_{\varepsilon}$, the metric can be written as one of the following three forms:
\begin{numcases}
{   \mathrm{d}{\sigma}^{2} \Big|_{U_{\varepsilon}} =   }   (\beta+1)^{2}|z|^{2\beta}|\mathrm{d}z|^2    &$\beta\in\mathbb{R}\setminus\{-1\}$, \\
c^2|z|^{-2}|\mathrm{d}z|^2                                         & $c>0$,     \\
{\left|\frac{\nu}{z}-\frac{n}{z^{n+1}}\right|^2|\mathrm{d}z|^2 }   &  $\nu>0$,  $n\in\mathbb{Z}$, $n\geq1$.
\end{numcases}
The constants $\beta, c, n, \nu$ above are all unique. Furthermore,  the coordinate $z$ is unique up to a rotation $z\mapsto\lambda z,  |\lambda|=1$ in the first form when $\beta\in((-\infty, -1)\setminus\mathbb{Z})\cup(-1, +\infty)$, and up to a scalar $z\mapsto\mathrm{p} z,  \mathrm{p}\in \mathbb{C}^*=\mathbb{C}\setminus\{0\} $ in the second form.
\end{theorem}

\setlength{\parskip}{1em}

\noindent \textbf{Remark 1.3.} We will see in Remark 2.1 that the three forms do not coincide with each other.

\setlength{\parskip}{1em}

\noindent \textbf{Remark 1.4.} The area growth condition is independent of the choice of the origin-preserving holomorphic coordinates. This is merely the change of variable formula in $\mathbb{R}^2$. So we can use the local expression of the metric in Theorem 1.2 to compute the area  growth rate corresponding to each form. It is easy to see that when $r$ goes to infinity,
\begin{equation}
\textrm{Area}(\Delta(0, \frac{1}{r}, \varepsilon)) =
\left\{ \begin{array}{ll}
O(1)  \ \textrm{when}\ \beta>-1, \ \textrm{and}\    O(r^{-2(\beta+1)})  \ \textrm{when}\ \beta<-1,  \\
O( \mathrm{ln}\; r), \\
O( r^{2n}).
\end{array} \right.  \notag
\end{equation}
So we immediately obtain a conclusion similar to Theorem 1.1:

\setlength{\parskip}{1em}

\noindent \textbf{Corollary 1.5.} \emph{   The isolated singularities of a flat metric with finite area on a Riemann surface must be conical singularities.  }

\setlength{\parskip}{1em}

\setlength\parindent{2em}The rest of this manuscript is organized as follows. In Section 2 we prove Theorem 1.2 using Complex Analysis. We give some complementary materials to the proof of Theorem 1.2 in the last section.



\clearpage

\setlength\parindent{0em}\textbf{\Large{2\ \ \  Proof of the main theorem}}

\setlength\parindent{2em} The rough idea of the proof is firstly to construct an  orientation-preserving local isometry from the punctured disk to the flat space form $\mathbb{C}$, then use the area condition to prove that the origin cannot be an essential singularity of a function appearing in the isometry, and finally choose suitable holomorphic coordinates to simplify the metric.

\setlength{\parskip}{1em}

\setlength\parindent{0em}\emph{Proof of theorem 1.2}. Let $L=\{y\in\mathbb{C}\mid Re\ {y}<0\}$ be the left half complex plane. Then $\mathrm{exp}:L\rightarrow{\Delta}^{*},\ y\mapsto{e}^y$ is a universal covering map of ${\Delta}^{*}$. Let $\mathrm{d}\tilde{\sigma}^2={\mathrm{exp}}^{*}(\mathrm{d}{\sigma}^{2})$, then $\mathrm{d}\tilde{\sigma}^2$ is a conformal metric on $L$ of Gauss curvature zero, and $\mathrm{d}\tilde{\sigma}^2$ is invariant under the deck transformation $y\mapsto{y+2{\pi}\mathrm{i}}$. Since $L$ is simply connected, it follows that there exists a well-defined holomorphic function $\xi:L\rightarrow\mathbb{C}$ so that
$\xi^*(\mathrm{d} \sigma_{0}^2)=\mathrm{d}\tilde{\sigma}^2$, where $\mathrm{d} \sigma_{0}^2$ is the Euclidean metric on the complex plane.
The invariance of $\mathrm{d}\tilde{\sigma}^2$ implies that there exists a point $(\theta,y_{0})\in[0,2\pi)\times\mathbb{C}$ such that
$\xi(y+2\pi\mathrm{i})=e^{\mathrm{i}\theta}\xi(y)+y_{0}$.
Denote $\theta=2\pi\alpha,\,0\le\alpha<1$.

\setlength{\parskip}{0.5em}

$\mathrm{CASE}\ 1. \ \ \theta\ne0, \ \mathrm{i.e.},\ 0<\alpha<1$. By composing a translation
$T:\mathbb{C}\rightarrow\mathbb{C}, \zeta\mapsto\zeta-\frac{y_{0}}{1-e^{\mathrm{i}\theta}}$
and replacing $\xi$ by $T\circ\xi$, we may assume that $\xi$ actually satisfies
$\xi(y+2\pi\mathrm{i})=e^{\mathrm{i}\theta}\xi(y)=e^{\mathrm{i}2\pi\alpha}\xi(y)$.
It follows that there exists a well-defined holomorphic function $\psi(\omega)$ on $\Delta^*$ such that
$\psi(e^y)=e^{-{\alpha}y}\xi(y)$.
It follows that there exists a multivalued locally univalent holomorphic function $f(\omega)=\omega^\alpha\psi(\omega)$ on $\Delta^*$ such that
$f\circ\mathrm{exp}=\xi$.
So we have
$f^*(\mathrm{d} \sigma_{0}^2)=\mathrm{d}{\sigma}^{2}$
($f$ is called the \emph{developing map} of $\mathrm{d}{\sigma}^{2}$).
It follows that
\begin{align}
\mathrm{d}{\sigma}^{2}
&=|\mathrm{d}f|^2=|\mathrm{d}(\omega^\alpha\psi(\omega))|^2              \notag \\
&=|\alpha{\omega}^{\alpha-1}\psi(\omega)+\omega^{\alpha}\psi'(\omega)|^2|\mathrm{d}\omega|^2       \notag\\
&=|\omega|^{2(\alpha-1)}|\alpha\psi(\omega)+\omega\psi'(\omega)|^2|\mathrm{d}\omega|^2,             \notag
\end{align}
where the last expression shows that the middle expressions are well defined even though $f(\omega)=\omega^\alpha\psi(\omega)$ is a multivalued function on $\Delta^*$. Next, we are going to use the hypothesis that the $\mathrm{d}{\sigma}^{2}$-area near $\omega=0$ satisfies the polynomial growth condition  to prove that $\psi$ is meromorphic at $\omega=0$. Equivalently, $\omega=0$ cannot be an essential singularity of $\psi$. Otherwise, the Laurent series
$\displaystyle \psi(\omega)=\sum_{n=-\infty}^{+\infty} a_n {\omega}^n$
would contains infinitely many nonzero coefficients $a_{-n}\ (n>0) $, which implies that the function
$\displaystyle \alpha\psi(\omega)+\omega\psi'(\omega)=\sum_{n=-\infty}^{+\infty} (\alpha+n)a_n {\omega}^n$ also has an essential singularity at
$\omega=0$, and its regular part has the same convergence radius (denote by $R$) as $\psi$'s. It follows that
\begin{equation}
\lim_{r \to +\infty}
\frac{\max \limits_{|\omega|=\frac{1}{r}}|\alpha\psi(\omega)+\omega\psi'(\omega)|}{r^{\eta}}
= +\infty, \ \ \forall \eta>0.             \notag
\end{equation}
Namely, $\forall{n}>0,\ n\in\mathbb{Z},\ \exists{r_{0}}>\frac{3}{R}$ such that
\begin{equation}
\max \limits_{|\omega|=\frac{1}{r}}|\alpha\psi(\omega)+\omega\psi'(\omega)|^2>3^{2-2\alpha}2^{n+2\alpha}r^{n+2\alpha}, \ \ \forall{r}>r_{0}.     \notag
\end{equation}
For ${r}>2r_{0}$, we consider the circle $\partial\Delta(0, \frac{2}{r})=\{|\omega|=\frac{2}{r}\}$ and suppose that $|\alpha\psi(\omega)+\omega\psi'(\omega)|^2$ gets its maximum value at $\omega=\omega_{0}\,(|\omega_{0}|=\frac{2}{r})$ on
$\partial\Delta(0, \frac{2}{r})$. Then on the disk $\Delta(\omega_{0}, \frac{1}{r})=\{|\omega-\omega_{0}|\leq\frac{1}{r}\}$, the mean value property yields
\begin{equation}
|\alpha\psi(\omega_{0})+\omega_{0}\psi'(\omega_{0})|^2\le\frac{1}{\pi(\frac{1}{r})^2}
\int\!\!\!\int_{\Delta(\omega_{0},\frac{1}{r})} |\alpha\psi(\omega)+\omega\psi'(\omega)|^2
\, \mathrm{d} x \mathrm{d} y.             \notag
\end{equation}
Combined with
\begin{equation}
|\alpha\psi(\omega_{0})+\omega_{0}\psi'(\omega_{0})|^2>3^{2-2\alpha}2^{n+2\alpha}(\frac{r}{2})^{n+2\alpha}=3^{2-2\alpha}r^{n+2\alpha},             \notag
\end{equation}
we have
\begin{align}
\mathrm{Area}(\Delta(0,\frac{1}{r},R))
&{\ge}  \mathrm{Area}(\Delta(\omega_{0},\frac{1}{r}))     \notag\\
&=\int\!\!\!\int_{\Delta(\omega_{0},\frac{1}{r})} |\omega|^{2(\alpha-1)}|\alpha\psi(\omega)+\omega\psi'(\omega)|^2\,\mathrm{d}x \mathrm{d}y    \notag\\
&>(\frac{3}{r})^{2(\alpha-1)}{\pi}(\frac{1}{r})^{2}3^{2-2\alpha}r^{n+2\alpha}             \notag\\
&={\pi}r^n,                          \notag
\end{align}
which contradicts the hypothesis in the theorem for the arbitrariness of $n$. So $\omega=0$ is either a pole or a removable singularity of $\psi$. It follows that there exists an integer $n$ and an $\varepsilon>0$ so that $\psi(\omega)=\omega^n e^{g(\omega)}$ on
$U_{\varepsilon}=\{ |\omega|<\varepsilon\}$, where $g(\omega)$ is holomorphic on $U_{\varepsilon}$. So $\mathrm{d}{\sigma}^{2} \Big|_{U_\varepsilon} =|\mathrm{d}(\omega^{\alpha+n}e^{g(\omega)})|^2$. We set $\beta=\alpha+n-1$ and
$z(\omega)={\omega}e^{\frac{g(\omega)}{\beta+1}}$, then we easily compute that
\begin{equation}
\mathrm{d}{\sigma}^{2} \Big|_{U_\varepsilon}=(\beta+1)^{2}|z|^{2\beta}|\mathrm{d}z|^2.              \notag
\end{equation}
Here $\beta\notin \mathbb{Z}$ for $0<\alpha<1$.

$\mathrm{CASE}\ 2. \ \ \theta=0, \ \mathrm{i.e.},\ \alpha=0$. So $\xi(y+2\pi\textrm{i})=\xi(y)+y_{0}$. Then there exists a well-defined holomorphic function
$\psi(\omega)$ on $\Delta^*$ such that $\psi(e^y)=\xi(y)-\frac{y_{0}}{2\pi\mathrm{i}}y$. It follows that the developing map
$f(\omega)=\frac{y_{0}}{2\pi\mathrm{i}}\ \mathrm{log}\;{\omega}+\psi(\omega)$ satisfies $f\circ\mathrm{exp}=\xi$. Denote $c=\frac{y_{0}}{2\pi\mathrm{i}}$, then
\begin{equation}
\mathrm{d}{\sigma}^{2}=|\mathrm{d}f|^2=\left| \frac{c}{\omega}+\psi'(\omega)  \right|^2|\mathrm{d} \omega|^2
=\left|  \frac{c+\omega\psi'(\omega)}{\omega}   \right|^2| \mathrm{d} \omega|^2.              \notag
\end{equation}

\setlength\parindent{2em}If $\omega=0$ is an essential singularity of $\psi$, then $\mathrm{Area}(\Delta(0,\frac{1}{r},R))$ cannot be controlled by any power of r. The discussion is analogous to $\mathrm{CASE}\ 1$. So $\omega=0$ is at worst a pole of $\psi$. Following the notations before, we have
$\psi(\omega)=\omega^n e^{g(\omega)}$ on $U_{\varepsilon}$.

\romannumeral1){} \ $y_{0}=0$, $\mathrm{i.e.},\  c=0$. Then $\mathrm{d}{\sigma}^{2}  \Big|_{U_\varepsilon} =|\mathrm{d}(\omega^{n}e^{g(\omega)})|^2$.

\textcircled{1}{} \ $n\ne0$. We set $\beta=n-1$ and $z(\omega)={\omega}e^{\frac{g(\omega)}{n}}$. Then
\begin{equation}
\mathrm{d}{\sigma}^{2} \Big|_{U_\varepsilon}=(\beta+1)^{2}|z|^{2\beta}|\mathrm{d}z|^2.   \notag
\end{equation}
Here $\beta\in \mathbb{Z}$ but $\beta\ne -1$.

\textcircled{2}{} \ $n=0$. Then $\omega=0$ is a removable singularity of $\psi$ and $\psi(0)\ne0$. So $\psi$ can be expanded into power series
$\displaystyle \psi(\omega)=\sum_{n=0}^{+\infty} a_n {\omega}^n, a_{0}\ne0$.
Since the Euclidean metric on the complex plane is translation-invariant, by composing a translation $\tilde{T}:\zeta\mapsto\zeta-a_{0}$, we have $\mathrm{d}{\sigma}^{2}\Big|_{U_\varepsilon} =|\mathrm{d}(\tilde{T}\circ\psi)|^2$. Now $\displaystyle (\tilde{T}\circ\psi)(\omega)=\omega^k \sum_{n=0}^{+\infty} a_{k+n}\; {\omega}^{n}$, where $a_{k}$ is the first nonzero coefficient after $a_{0}$, the discussion goes back to \textcircled{1}.

\setlength\parindent{2em}The constant $\beta$ in $(1)$ is unique since it characterizes the growth rate of $\mathrm{d}{\sigma}^{2}$-area near the origin. In fact, let $\tilde{z}$ be another origin-preserving holomorphic coordinate on $U_{\varepsilon}$ such that
\begin{equation}
\mathrm{d}{\sigma}^{2} \Big|_{U_\varepsilon}=(\tilde{\beta}+1)^{2}|\tilde{z}|^{2\tilde{\beta}}|\mathrm{d}\tilde{z}|^2, \tilde{\beta}\neq-1.
\end{equation}
Then there exists a non-vanishing holomorphic function $h(z)$ near the origin such that $\tilde{z}=zh(z)$(we will frequently use this argument later). It follows that on $U_\varepsilon$(possibly smaller),
\begin{align}
\mathrm{d}{\sigma}^{2} \Big|_{U_\varepsilon}
&=(\tilde{\beta}+1)^{2}|\tilde{z}|^{2\tilde{\beta}}|\mathrm{d}\tilde{z}|^2    \notag \\
&=(\tilde{\beta}+1)^{2}\left|z h(z)\right|^{2\tilde{\beta}}    \left| (zh'(z)+h(z))\right|^2|\mathrm{d}z|^2   \notag \\
&=(\tilde{\beta}+1)^{2}\left|h(z)\right|^{2\tilde{\beta}} \left| (zh'(z)+h(z))\right|^2
\left|z \right|^{2\tilde{\beta}}|\mathrm{d}z|^2.
\end{align}
We compare $(1)$ and $(5)$, and notice that $(\tilde{\beta}+1)^{2}\left|h(z)\right|^{2\tilde{\beta}} \left| (zh'(z)+h(z))\right|^2$ is non-vanishing on $U_{\varepsilon}$ if $\varepsilon$ is small enough, therefore controlled by two positive numbers, which implies $\tilde{\beta}=\beta$.

\setlength{\parskip}{0em}

It is easy to see from $(1)$ and $(4)$ that $f_{1}(z)=z^{\beta+1}$ and $f_{2}(\tilde{z})=\tilde{z}^{\tilde{\beta}+1}=\tilde{z}^{\beta+1}$ are two developing maps of $\mathrm{d}{\sigma}^{2}\Big|_{U_\varepsilon}$ under the coordinate $z$ and $\tilde{z}$ respectively. So there exist $\zeta_0, \lambda_1\in \mathbb{C}, |\lambda_1|=1$ such that $f_{2}=\lambda_1 f_{1} + \zeta_0$ on $U_{\varepsilon}$. Substituting $\tilde{z}=zh(z)$ we get
\begin{equation}
(h^{\beta+1}-\lambda_1)z^{\beta+1}=\zeta_0.
\end{equation}

If $\beta>-1$, then $h^{\beta+1}-\lambda_1=0$, so $h=\lambda, |\lambda|=1$; If $\beta<-1$ and $\beta\notin \mathbb{Z}$ we have
$(h^{\beta+1}-\lambda_1)=\zeta_0 z^{-(\beta+1)}$, where the left hand side has single-valued branches near the origin since $h(0)\neq 0$ while the right hand side does not unless $\zeta_0 =0$, which still implies $h=\lambda, |\lambda|=1$; If $\beta<-1$ is a negative integer, then we directly solve the equation $(6)$ to obtain
\begin{equation}
h(z)=(\zeta_0 z^{-(\beta+1)}+\lambda_1)^{\frac{1}{\beta+1}}.
\end{equation}

\setlength{\parskip}{1em}

\romannumeral2){} \ $y_{0}\ne0$, then\ $c=\frac{y_{0}}{2\pi\mathrm{i}}\ne0$.

\setlength{\parskip}{0.5em}

\textcircled{1}{} \ $n\ge0$, then $\omega=0$ is a removable singularity of $\psi$. Setting
$z={\omega}e^{\frac{\psi(\omega)}{c}}$ to get
\begin{equation}
\mathrm{d}{\sigma}^{2} \Big|_{U_\varepsilon}=\left| \frac{c}{z} \right|^{2} |\mathrm{d}z|^2,     \notag
\end{equation}
which is the second form $(2)$ if we still denote $|c|>0$ by $c$.

\setlength{\parskip}{0em}

It is obvious that a scalar  $z\mapsto \mathrm{p}z, \mathrm{p}\in \mathbb{C}^*$ keeps the second form. We want to show the inverse is also true. Let $\tilde{z}=zh(z)(h$ is non-vanishing and holomorphic near the origin) be another coordinate such that
\begin{equation}
\mathrm{d}{\sigma}^{2} \Big|_{U_\varepsilon}=\tilde{c}^2  \left| \tilde{z} \right|^{-2} |\mathrm{d}\tilde{z}|^2, \ \ \tilde{c}>0,            \notag
\end{equation}
Then $f_{1}(z)=c\ \mathrm{log}\; z$ and $f_{2}(\tilde{z})=\tilde{c}\ \mathrm{log}\; \tilde{z}$ are two developing maps of $\mathrm{d}{\sigma}^{2} \Big|_{U_\varepsilon}$ under the coordinate $z$ and $\tilde{z}$ respectively. So there exist $\zeta_0, \lambda\in\mathbb{C},|\lambda|=1$ such that
 $f_{2}=\lambda f_{1} + \zeta_0$ on $U_{\varepsilon}$. Substituting $\tilde{z}=zh(z)$ we get
\begin{equation}
(\tilde{c}-\lambda  c) \ \mathrm{log}\; z=-\tilde{c} \ \mathrm{log}\; h +\zeta_0
\end{equation}
The right hand side of $(8)$ has single-valued branches near the origin while the left hand side does not unless $\tilde{c}-\lambda  c=0$, which implies $\lambda=1, \tilde{c}=c$. Substituting this in $(8)$ to get $c \ \mathrm{log}\; h=\zeta_0$, so $h$ is a constant.

\setlength{\parskip}{0.5em}

\textcircled{2}{} \ $n<0$, then $\omega=0$ is a pole of $\psi$. This case echoes the study of Strebel in [3, Chapter\uppercase\expandafter{\romannumeral3}, \S 6, 6.4]. For convenience, we replace $n$ by $-n$ but still denote by $n>0$. Expand $\psi$ into Laurent series $\psi(\omega)=\frac{1}{\omega^{n}} \left( c_0+c_1 \omega +c_2 \omega^2+\cdots  \right), c_0\neq0.$
Then near the origin the developing map
\begin{equation}
f(\omega)=c \ \mathrm{log}\;\omega+\frac{1}{\omega^{n}} \left( c_0+c_1 \omega +c_2 \omega^2+\cdots  \right).
\end{equation}
Denote $c=\nu e^{\mathrm{i}\vartheta}, \nu>0$ and $g(\omega)=c_0+c_1\omega+c_2\omega^2+\cdots$.
Naturally, we want to find a local holomorphic coordinate $ z(\omega)=\omega e^{\phi(\omega)}$, where $\phi$ is an undetermined holomorphic function near the origin, so that
\begin{equation}
f(z)=c \ \mathrm{log}\;z+\frac{e^{\mathrm{i}\vartheta}}{z^{n}} +\eta , \ \ \eta\in\mathbb{C}.
\end{equation}
Then
\begin{equation}
\mathrm{d}{\sigma}^{2} \Big|_{U_{\varepsilon}}=|\mathrm{d}f|^2=\left|  \frac{\nu}{z}-\frac{n}{z^{n+1}}  \right|^2  |\mathrm{d} z   |^2,
\end{equation}
which is the third form  in $(3)$.

\setlength{\parskip}{0em}

Substituting $ z=\omega e^{\phi(\omega)}$ in $(10)$ and then comparing $(10)$ and $(9)$, we obtain
\begin{equation}
c\ \mathrm{log}(\omega e^{\phi(\omega)})+e^{\mathrm{i}\vartheta}(\omega e^{\phi(\omega)})^{-n}+\eta=c\ \mathrm{log}\;\omega+\omega^{-n}g(\omega).
\end{equation}
Reorganise $(12)$,
\begin{equation}
g(\omega)=e^{\mathrm{i}\vartheta} e^{-n\phi}+\omega^n(c\phi+\eta).
\end{equation}
From $(13)$ we know that $\phi(0)$ satisfies $e^{\mathrm{i}\vartheta} e^{-n\phi(0)}=g(0)=c_0$.
Taking derivative of $(13)$ yields
\begin{equation}
\phi'=\frac{g'-n\omega^{n-1}(c\phi+\eta)}{c\omega^n-ne^{\mathrm{i}\vartheta}e^{-n\phi}}.
\end{equation}
By Cauchy-Kovalevskaya theorem, $(14)$ has a local solution $\phi(\omega)$ for any initial value $\phi(0)$. By integration we get $(13)$, up to an additive constant, which must be zero since $e^{\mathrm{i}\vartheta} e^{-n\phi(0)}=g(0)$.

It is easy to see from $(3)$ that $n+1$ determines the growth rate of $\mathrm{d}{\sigma}^{2}$-area near the origin. So $n$ is unique, of which the discussion is analogous to that of $\beta$ in $(1)$. It remains to show the uniqueness of the positive number $\nu=|c|=\frac{|y_0|}{2\pi}$. Similar to the previous discussion, let $\tilde{z}=zh(z), h(0)\neq0$ be another origin-preserving biholomorphic coordinate such that
\begin{equation}
\mathrm{d}{\sigma}^{2}\Big|_{U_\varepsilon} = \left|  \frac{\tilde{\nu}}{\tilde{z}} -\frac{n}{\tilde{z}^{n+1}} \right|^2  |\mathrm{d}\tilde{z}|^2,  \tilde{\nu}>0.                    \notag
\end{equation}
Obviously $f_1(z)=\nu\ \mathrm{log}\; z+z^{-n}$ and $f_2(\tilde{z})=\tilde{\nu}\ \mathrm{log}\; \tilde{z}+\tilde{z}^{-n}$ are two developing maps of
$\mathrm{d}{\sigma}^{2} \Big|_{U_\varepsilon}$ under the coordinate $z$ and $\tilde{z}$ respectively. So there exist $\zeta_0, \lambda\in\mathbb{C}, |\lambda|=1$ such that $f_2(\tilde{z})=\lambda  f_1(z)+\zeta_0$. Substituting $\tilde{z}=zh(z)$ we have
\begin{equation}
\tilde{\nu}\ \mathrm{log}\; z+ \tilde{\nu}\ \mathrm{log}\; h+ (zh)^{-n}=\lambda\nu \ \mathrm{log}\; z+\lambda  z^{-n}+\zeta_0.
\end{equation}
Since $\mathrm{log}\; z$ is multi-valued near the origin, then $\nu=\lambda\tilde{\nu}$, which implies  $\lambda=1, \nu=\tilde{\nu}$. Thus the uniqueness of $\nu$ has been proved. Thus we complete the proof of the whole theorem. \qed

\setlength{\parskip}{1em}

\noindent \textbf{Remark 2.1.} Readers may want to try making the third form simpler, just like the previous two forms only containing one term in the absolute value symbol. But the efforts would be in vain because the developing maps of the three forms cannot differ from each other by a rotation and a translation. The discussion is analogous to that of the uniqueness of the constant $c, \nu$ and the coordinate transformation, which frequently appears in the proof.

\setlength{\parskip}{1em}

\noindent \textbf{Remark 2.2.} In $\mathrm{CASE}\ 2$, $|y_0|$ is invariant, which is an implication of the uniqueness of $c$ and $\nu$ (both equal $\frac{|y_0|}{2\pi}$). In fact, let $\tilde{\xi}$ be another holomorphic function from $L$ to $\mathbb{C}$ such that $\tilde{\xi}^*(\mathrm{d} \sigma_{0}^2)=\mathrm{d}\tilde{\sigma}^2$, then there exist $\zeta_0, \lambda\in\mathbb{C}, |\lambda|=1$ such that
$\tilde{\xi}(y)=\lambda\xi(y)+\zeta_0$. Suppose  $\tilde{\xi}(y+2\pi\mathrm{i})=e^{\mathrm{i}\tilde{\theta}}\tilde{\xi}(y)+\tilde{y}_0$, then we have
\begin{equation}
\lambda\xi(y+2\pi\mathrm{i})+\zeta_0=e^{\mathrm{i}\tilde{\theta}}(\lambda\xi(y)+\zeta_0)+\tilde{y}_0.
\end{equation}
Substituting $\xi(y+2\pi\mathrm{i})=\xi(y)+y_{0}$ in $(16)$ we get $e^{\mathrm{i}\tilde{\theta}}=1, \tilde{y}_0=\lambda  y_{0}$, so $|\tilde{y}_0|=|y_0|$.

\noindent \textbf{Remark 2.3.} Under the coordinate transformation $ z=e^{-\rho}e^{\mathrm{i}\theta} $, The second form
\begin{align}
\mathrm{d}{\sigma}^{2} \Big|_{U_\varepsilon} &=c^2|z|^{-2}|\mathrm{d}z|^2=c^2  \left|  e^{-\rho} \right|^{-2}
\left|  -e^{-\rho}e^{\mathrm{i}\theta} \mathrm{d}\rho + e^{-\rho}e^{\mathrm{i}\theta}\mathrm{i} \mathrm{d}\theta \right|^2      \notag \\
&=c^2 e^{2\rho}e^{-2\rho} \left| -\mathrm{d}\rho+ \mathrm{i} \mathrm{d}\theta  \right|^2 =c^2(\mathrm{d}\rho^2+\mathrm{d}\theta^2).               \notag
\end{align}
This is the Euclidean metric on a cylinder.

\noindent \textbf{Remark 2.4.} The coordinate transformation of the third form is not unique because $(11)$ holds for any value of $\eta$ in $(10)$. However, from another perspective, we still can say something about the `uniqueness'. We fix the coordinate $z$ and denote by $\mathcal{M}_i, i=1, 2, 3$ the set of all origin-preserving holomorphic coordinates keeping the metric form in $(i)$. Then the main theorem claims that
$\mathcal{M}_1\simeq S^1$  when $\beta>-1$ or
$\beta<-1, \beta\notin\mathbb{Z}$ and $\mathcal{M}_2\simeq \mathbb{C}^*$. The multiplication structures on $S^1$ and $\mathbb{C}^*$ correspond to the composition of the coordinates in $\mathcal{M}_1$ and $\mathcal{M}_2$ respectively. From $(7)$ we know that there is a bijection from $\mathcal{M}_1$ to $S^1\times\mathbb{C}$ when $\beta<-1$ and $\beta\in\mathbb{Z}$. The composition of coordinates endows $S^1\times\mathbb{C}$ a (non-commutative) multiplication structure:
$(\lambda_1, \zeta_1)\times(\lambda_2, \zeta_2)=(\lambda_1\lambda_2, \lambda_2\zeta_1+\zeta_2)$.
It is natural to ask about $\mathcal{M}_3$. We will show that locally $\mathcal{M}_3$ looks like $\mathbb{C}$, which is not surprising because of the one degree of freedom $\eta$.

\setlength{\parskip}{0em}

We need to clarify how $\tilde{z}$ (equivalently, $h(z)$) depends on $z$. So we go back to $(15)$. We have proved that $\lambda=1, \nu=\tilde{\nu}$ in (15), which yields
\begin{equation}
\nu\ \mathrm{log}\; h+ (zh)^{-n}=  z^{-n}+\zeta_0.
\end{equation}
As $z$ approaches to zero we get $h(0)^{-n}=1$, hence $h(0)=e_k:=e^{\frac{2k\pi\mathrm{i}}{n}}, k\in\{0,1,2, \cdots, n-1\}$.
Taking derivative of $(17)$ and reorganizing to obtain
\begin{numcases}
{      }  \begin{split}
&h'=\frac{nh(1-h^n)}{z(\nu z^n h^n-n)}, \\
&h(0)=e_k.
\end{split}
\end{numcases}
Let $\tilde{h}=h-e_k$, then $(18)$ yields
\begin{numcases}
{      }  \begin{split}
&z\tilde{h}'=\tilde{F}(z,\tilde{h}), \\
&\tilde{h}(0)=0,
\end{split}
\end{numcases}
where $\tilde{F}(z,\tilde{h})=\frac{n(\tilde{h}+e_k) (1-(\tilde{h}+e_k)^n)}{\nu z^n (\tilde{h}+e_k)^n-n}$ is holomorphic at $(0,0)$. Expand $\tilde{F}$ into power series
$$
\tilde{F}(z,\tilde{h})=n\tilde{h}+\sum_{\mbox{\tiny$\begin{array}{c}
j+l\geq2\\
\mathrm {if}\ l=0 \ \mathrm {then} \ j\geq2n  \end{array}$}}c_{jl}z^j \tilde{h}^l.
$$

If $n=1$, by theorem $11.1.4$ in [4], there exist infinitely many local solutions $\tilde{h}$ for the derivative $\tilde{h}'(0)$ is arbitrary. i.e., if we expand $\tilde{h}$ into power series $\tilde{h}(z)=a_1 z + a_2 z^2 +\cdots$, then $a_1$ is arbitrary and once $a_1$ has been chosen, $a_j\textrm{'}s(j\geq2)$ are all determined by $a_1$ and ${c_{jl}}\textrm{'}s$. And we get $h=\tilde{h}+e_k=\tilde{h}+1=1+a_1z+a_2 z^2+\cdots$.

If $n\geq2$, taking $\tilde{h}=z^{n-1}\hat{h}$ and substituting it into $(19)$, we get $\displaystyle  z\hat{h}'=\hat{h}+\sum_{j+l\geq2} \hat{c}_{jl}z^j \hat{h}^l$,
the same type of ODE as the one of $n=1$. So we have infinitely many solutions $\hat{h}=\hat{a}_1 z + \hat{a}_2 z^2+\cdots$, where $\hat{a}_1$ is arbitrary. Then
$h=\tilde{h}+e_k=z^{n-1}\hat{h}+e_k=e_k+\hat{a}_1z^n+\hat{a}_2 z^{n+1}+\hat{a}_3 z^{n+2}+\cdots$.

Summing up the above two cases, we obtain a bijection from $\mathcal{M}_3$ to $\{e_0, \cdots, e_{n-1}\}\times\mathbb{C}\simeq\mathop{\bigsqcup}\limits_n \mathbb{C}$ via the map
\begin{align}
\mathcal{M}_3=&\{  z\mapsto \tilde{z}=zh(z) \mid \   \tilde{z} \ \textrm{keeps the third form}  \}
\rightarrow  \{e_0, \cdots, e_{n-1}\}\times\mathbb{C}        \notag  \\
&h  \mapsto  (h(0), \frac{h^n(0)}{n!}).      \notag
\end{align}
The multiplication structure is given by $(e_k, a_n)\times(e_j, \tilde{a}_n)=(e_k e_j, e_k \tilde{a}_n+e_j a_n)$.\\

\setlength\parindent{0em}\textbf{\Large{3\ \ \  Some complementary materials}}

\setlength{\parskip}{1em}

3.1. {\it The area growth condition is independent of the choice of the origin-preserving holomorphic coordinates. i.e.,
let $z$ be another origin-preserving holomorphic coordinate near the origin, then we can find new $0<R_1<1, M_1>0$ such that}
$$\textrm{Area}\bigg(\Delta_z\Big(0, \frac{1}{r}, R_1\Big)\bigg)\leq M_1 r^N, \forall r>\frac{1}{R_1}.$$

\setlength\parindent{2em}We write $\omega=zh(z)$, where $h(z)$ is non-vanishing holomorphic in some closed neighborhood of $z=0$, say $\overline{U}_{R_2}=\{|z|\leq  R_2(<1)\}$, therefore $|h|$ has positive lower and upper bounds on $\overline{U}_{R_2}$, denoted by $B, A>0$ respectively. Then we have $0<B\leq\left| \frac{\omega}{z} \right|=|h|\leq  A$ on $\overline{U}_{R_2}$. Suppose
\begin{equation}
\mathrm{d}{\sigma}^{2} \Big|_{U_{R_2}} =e^{2\psi(z)} |\mathrm{d}z|^2 = e^{2\varphi(\omega)}|\mathrm{d}\omega|^2,  \notag
\end{equation}
then $e^{2\varphi(\omega)}=e^{2\psi(z(\omega))} |z'(\omega)|^2$.
Choose $R_1<\min\{ R_2, \frac{R}{A} \}$ and consider $\Delta_z(0, \frac{1}{r}, R_1)=\{ \frac{1}{r}<|z|<R_1 \}$ with $r>\frac{1}{R_1}$. Notice that
$\frac{1}{r}<|z|<R_1 \Rightarrow \frac{B}{r}<|\omega|<AR_1 $, which implies $\omega(\{ \frac{1}{r}<|z|<R_1 \})\subset\{ \frac{B}{r}<|\omega|<AR_1 \}$, we can see that
\begin{align}
\textrm{Area}(\Delta_z(0, \frac{1}{r}, R_1))&=
\int\!\!\!\int_{  \{ \frac{1}{r}<|z|<R_1 \}  } \frac{\mathrm{i}}{2} e^{2\psi(z)}
\, \mathrm{d} z \wedge \mathrm{d} \bar{z}      \notag\\
&=\int\!\!\!\int_{ \omega(\{ \frac{1}{r}<|z|<R_1 \})  } \frac{\mathrm{i}}{2} e^{2\psi(z(\omega))} |z'(\omega)|^2
\, \mathrm{d} \omega \wedge \mathrm{d} \bar{\omega}    \notag\\
&=\int\!\!\!\int_{ \omega(\{ \frac{1}{r}<|z|<R_1 \})  } \frac{\mathrm{i}}{2} e^{2\varphi(\omega)}
\, \mathrm{d} \omega \wedge \mathrm{d} \bar{\omega}     \notag \\
&\leq  \int\!\!\!\int_{ \{ \frac{B}{r}<|\omega|<AR_1 \}  } \frac{\mathrm{i}}{2} e^{2\varphi(\omega)}
\, \mathrm{d} \omega \wedge \mathrm{d} \bar{\omega}      \notag \\
&= \textrm{Area}(\Delta_\omega(0, \frac{1}{\frac{r}{B}}, AR_1(<R))) \leq M(\frac{r}{B})^N = \frac{M}{B^N}r^N.  \notag
\end{align}
So $N$ is invariant.

\noindent 3.2. {\it Explanation about Remark 2.1}

\setlength\parindent{2em} Assume that $\tilde{z}=zh(z), h(0)\neq0$ and $\zeta_0, \lambda\in\mathbb{C}, |\lambda|=1.$
If $z^{\beta+1}=\lambda c\ \mathrm{log}\; \tilde{z} + \zeta_0 =\lambda c\ \mathrm{log}\; z + \lambda c\ \mathrm{log}\; h+ \zeta_0,$  when $z=x>0$ small enough, the right hand side $\sim$ $\mathrm{ln}\; x$, which implies that $\beta<-1$ and so the left hand side $\sim$ $x^{\beta+1}$. Impossible! So (1) does not coincide with (2);
If $\nu\ \mathrm{log}\; z+z^{-n}=\lambda \tilde{z}^{\beta+1} +\zeta_0=\lambda h^{\beta+1}z^{\beta+1}+\zeta_0,$ the same argument shows that $\beta=-n-1$ and so $\nu\ \mathrm{log}\; z=\lambda h^{-n}z^{-n}-z^{-n}+\zeta_0.$  Impossible! So (1) does not coincide with (3);
If $\nu\ \mathrm{log}\; z+z^{-n}=\lambda c\ \mathrm{log}\; \tilde{z} + \zeta_0=\lambda c\ \mathrm{log}\; z + \lambda c\ \mathrm{log}\; h+ \zeta_0,$
the analogous discussion  implies that (2) does not coincide with (3).

\noindent 3.3. {\it Why is $n\tilde{h}$  the first term of the power series expansion of $\tilde{F}(z,\tilde{h})=\frac{n(\tilde{h}+e_k) (1-(\tilde{h}+e_k)^n)}{\nu z^n (\tilde{h}+e_k)^n-n}$
in Remark 2.4  ?}

\setlength\parindent{2em}We only need to compute the formal power series. Noticing that when $   (z,\tilde{h})\rightarrow(0, 0), (z,h)\rightarrow(0, e_k),  \frac{\nu}{n}z^nh^n =o(1), $ we have
\begin{equation}
\begin{aligned}
\frac{nh (1-h^n)}{\nu z^n h^n-n}&=\frac{-h (1-h^n)}{1-\frac{\nu}{n}z^nh^n}  \notag\\
&= (-h+h^{n+1})   \left(   1+ \frac{\nu}{n}z^nh^n +  (\frac{\nu}{n}z^nh^n )^2 +  (\frac{\nu}{n}z^nh^n )^3 + \cdots \right)   \notag\\
&=-h- \frac{\nu}{n}z^nh^{n+1} - \frac{\nu^2}{n^2}z^{2n}h^{2n+1} - \cdots    \notag\\
&\ \ \ +h^{n+1}+ \frac{\nu}{n}z^nh^{2n+1} + \frac{\nu^2}{n^2}z^{2n}h^{3n+1} +\cdots.
\end{aligned}
\end{equation}
Computing the two lowest-order terms  $-h+h^{n+1}$ to get
\begin{equation}
\begin{aligned}
\tilde{F}(z,\tilde{h})
&=-(\tilde{h}+e_k)- \frac{\nu}{n}z^n(\tilde{h}+e_k)^{n+1} - \frac{\nu^2}{n^2}z^{2n}(\tilde{h}+e_k)^{2n+1} - \cdots   \notag\\
&\ \ \ +(\tilde{h}+e_k)^{n+1}+ \frac{\nu}{n}z^n(\tilde{h}+e_k)^{2n+1} + \frac{\nu^2}{n^2}z^{2n}(\tilde{h}+e_k)^{3n+1} +\cdots    \notag\\
&=-\tilde{h}-e_k+\tilde{h}^{n+1}+(n+1)\tilde{h}^n e_k + \cdots + (n+1)\tilde{h} e_k^n + e_k^{n+1}     \notag\\
&\ \ \ -\frac{\nu}{n}z^n(\tilde{h}+e_k)^{n+1} - \frac{\nu^2}{n^2}z^{2n}(\tilde{h}+e_k)^{2n+1} - \cdots      \notag\\
&\ \ \ +\frac{\nu}{n}z^n(\tilde{h}+e_k)^{2n+1} + \frac{\nu^2}{n^2}z^{2n}(\tilde{h}+e_k)^{3n+1} +\cdots     \notag\\
&=n\tilde{h}+\sum\limits_{j+l\geq2} c_{jl} z^j \tilde{h}^l
\end{aligned}
\end{equation}

We notice that $-e_k\frac{\nu}{n} z^n$ in $- \frac{\nu}{n}z^n(\tilde{h}+e_k)^{n+1}$ balances $e_k\frac{\nu}{n} z^n$ in $+ \frac{\nu}{n}z^n(\tilde{h}+e_k)^{2n+1}$, which yields that if $l=0$ then $j\geq2n$. This is important to guarantee that when $n\geq2$, the transformation $\tilde{h}=z^{n-1}\hat{h}$ does not damage the holomorphicity of the relevant functions in the new ODE.

\noindent 3.4. {\it The multiplication induced by the composition of coordinates in remark 2.4}

\setlength\parindent{2em}When $\beta<-1$ and $\beta\in\mathbb{Z}$ in the first form,
\begin{equation}
\begin{aligned}
\tilde{z}&=zh_1(z)=z(\zeta_1 z^{-(\beta+1)}+\lambda_1)^{\frac{1}{\beta+1}},   \notag\\
\hat{z}
&=\tilde{z}h_2(\tilde{z})=\tilde{z}(\zeta_2 \tilde{z}^{-(\beta+1)}+\lambda_2)^{\frac{1}{\beta+1}}   \notag\\
&=z(\zeta_1 z^{-(\beta+1)}+\lambda_1)^{\frac{1}{\beta+1}}  \left( \zeta_2 z^{-(\beta+1)}(\zeta_1 z^{-(\beta+1)}+\lambda_1)^{-1}+\lambda_2  \right)^{\frac{1}{\beta+1}}  \notag\\
&=z \left(  \zeta_2 z^{-(\beta+1)}+\lambda_2(\zeta_1 z^{-(\beta+1)}+\lambda_1)   \right)^{\frac{1}{\beta+1}}   \notag\\
&=z  \left(  (\zeta_2+\lambda_2 \zeta_1) z^{-(\beta+1)}+\lambda_1\lambda_2   \right)^{\frac{1}{\beta+1}}.
\end{aligned}
\end{equation}

In the third form,
\begin{equation}
\begin{aligned}
\tilde{z}&=zh_1(z)=z(e_k+a_n z^n+a_{n+1} z^{n+1}+\cdots),      \notag\\
\hat{z}
&=\tilde{z}h_2(\tilde{z})=\tilde{z}(e_j+\tilde{a}_n \tilde{z}^n+\tilde{a}_{n+1} \tilde{z}^{n+1}+\cdots)   \notag\\
&=z(e_k+a_n z^n+a_{n+1} z^{n+1}+\cdots)\cdot             \notag\\
&\ \ \ \left(   e_j+\tilde{a}_n z^n(e_k+a_n z^n+\cdots)^n + \tilde{a}_{n+1} z^{n+1}(e_k+a_n z^n+\cdots)^{n+1} +\cdots \right)     \notag\\
&=z(e_k+a_n z^n+a_{n+1} z^{n+1}+\cdots)  \left(  e_j+\tilde{a}_n z^n(e_k^n+A(z)z^n) +B(z)z^{n+1}   \right)    \notag\\
&=z(e_k+a_n z^n+a_{n+1} z^{n+1}+\cdots) (e_j+\tilde{a}_n z^n+C(z)z^{n+1})      \notag\\
&=z \left(   e_k e_j+(e_k\tilde{a}_n+e_j a_n)z^n+\textrm{high-order terms}    \right).
\end{aligned}
\end{equation}

\clearpage

$$\textbf{Acknowledgements}$$

\noindent Xu would like to express his sincere gratitude to Professor Song Sun for the stimulating conversations in the summer of 2016, which motivated this manuscript. Xu is supported in part by the National Natural Science Foundation of China (No. 11571330) and the Fundamental Research Funds for the Central Universities. Li is support in part by the Research Training Group 1821 ``Cohomological Methods in Geometry".

{\sc \noindent Jin Li\\
Mathematical Institute, Albert Ludwigs University of Freiburg\\
Ernst-Zermelo-Str. 1\\
79104 Freiburg im Breisgau, Germany\\
jin.li@math.uni-freiburg.de
}\\

{\sc \noindent Bin Xu\\
Wu Wen-Tsun Key Laboratory of Math, USTC, CAS\\
School of Mathematical Sciences\\
University of Science and Technology of China\\
Hefei 230026 China\\
bxu@ustc.edu.cn}

\end{document}